\documentclass[12pt, oneside]{amsart}
\usepackage{epsfig, amssymb, amsfonts}
\usepackage[left=1.5in, right=1in, top=1in, bottom=1in]{geometry}

\newcommand{\skipline}{\vspace{12pt}}
\newtheorem{theorem}{Theorem}
\newtheorem{definition}{Definition}
\newtheorem{lemma}{Lemma}

\begin{document}

\pagenumbering{roman}
\pagestyle{plain}

\title{Time-Dependent Solutions of a Discrete Schr\"{o}dinger's Equation}
\author{Nigie Shi}
\date{\today}

\begin{center}

\begin{Large}Time-Dependent Solutions of a Discrete Schr\"{o}dinger's Equation\end{Large}\\
\skipline
By\\
\skipline
NIGIE SHI\\
B.S. (University of Wyoming) 2001\\
\skipline
THESIS\\
\skipline
Submitted in partial satisfaction of the requirements for the degree of\\
\skipline
MASTER OF SCIENCE\\
\skipline
in\\
\skipline
APPLIED MATHEMATICS\\
\skipline
in the\\
\skipline
OFFICE OF GRADUATE STUDIES\\
\skipline
of the\\
\skipline
UNIVERSITY OF CALIFORNIA,\\
\skipline
DAVIS\\
\vspace{1in}
Approved:\\
\skipline
\rule{2.5in}{1pt}\\
\skipline
\rule{2.5in}{1pt}\\
\skipline
\rule{2.5in}{1pt}\\
\skipline
\rule{2.5in}{1pt}\\
\skipline
Committee in Charge\\
\skipline
2004\\

\end{center}

\renewcommand{\baselinestretch}{1.6}\small\normalsize

\newpage
\tableofcontents

%

\newpage
\section*{Acknowledgments}

I would like to thank for Professor Bruno Nachtergaele for giving
me this interesting project and being my thesis advisor and
Professor Albert C. Fannjiang and Professor Alexander Soshnikov
who are working in related fields for kindly agreeing to be
members of my thesis committee and review my thesis. I especially
like to show my great appreciation to Professor Bruno
Nachtergaele's patient guidance for helping me complete this
project. The process is very invaluable. Professor Albert C.
Fannjiang and Professor Alexander Soshnikov are also my
instructors for the two graduate preliminary math courses MAT 119
and MAT 203, respectively. I have learned a lot from them during
my first year of graduate study at Uuniversity of California,
Davis. Also, I would like to thank for Arthur Cheng, Jeremy Clark,
Ben-Shan Liao, and Yuan-Kai Huang who are my colleagues for
consulting related materials of my project in both analysis and
numerics.

Moreover, mathematics is always my primary interest including
researching and teaching that I would like to work on both of them
for my future career and I sincerely appreciate Professor Bruno
Nachtergaele, Professor Alex Mogilner, Professor Elbridge Gerry
Puckett, Doctor Duane Kouba, and Professor Jim Diederich for
writing my letters of recommendation.

Finally, I would like to thank for my parents Mei-Pin Shi and
Ling-Ju Lin of supporting me to come to United States of America
for my college education at University of Wyoming. Based on my
hard work during undergraduate and their constant encouragement, I
am able to go to University of California, Davis for my graduate
study.

\newpage
\pagestyle{headings}
\pagenumbering{arabic}

\section{Introduction}

A very useful method of studying linear operators by decomposing
the space on which they act into invariant subspaces is known as
spectral theory. An example of an application of spectral theory
is the problem of finding a set of eigenvectors or diagonalizing a
linear map on an infinite-dimensional space [HN]. When a finite
dimensional linear operator is diagonalized, there exist a set of
eigenvalues and their corresponding eigenvectors. Along the
directions of an eigenvector with its given eigenvalue, the action
of the operator is just multiplication by the eigenvalue. The
spectrum contains the set of eigenvalues that is also called point
spectrum. In infinite dimensional case, the structure of the
spectrum will often be more complicated such that there may exist
a continuous spectrum or residual spectrum which do not contain a
set of eigenvalues and are different from the point spectrum. In
my thesis, I study the spectrum of an operator based on some
numerical results of it and also properties of the discrete
Laplacian operator.

The discrete Laplacian operator $\Delta$ is defined on
$\ell^{2}(\mathbb{Z})$ by
\begin{equation*}
(\Delta x)_{k}=(x_{k-1}+x_{k+1})-2x_{k},
\end{equation*}
where $x=(x_k)_{k=-\infty}^{\infty} \in \ell^{2}(\mathbb{Z})$.
Moreover, three features of $\Delta$ are described in a theorem as
follows:
\begin{theorem}
Let $\mathcal{S}$ be the right shift operator and
$\mathcal{S^{*}}$ be the left shift operator.
$\Delta=\mathcal{S}+\mathcal{S^{*}}-2\mathcal{I}$. The spectrum of
$\Delta$ is entirely continuous and consists of the interval
$[-4,0]$.
\end{theorem}

\textbf{Proof.} First, $\|\Delta\| \leq \|\mathcal{S}\| +
\|\mathcal{S^{*}}\| + 2\|\mathcal{I}\| = 4$ where $\|\cdot\|$
denotes the operator norm. Second, $\Delta =
-(\mathcal{S}-\mathcal{I})(\mathcal{S^{*}}-\mathcal{I})=-(\mathcal{S}-\mathcal{I})(\mathcal{S}-\mathcal{I})^{*}
\leq 0$. These two arguments imply that $\sigma(\Delta) \subseteq
[-4,0]$ where $\sigma(\Delta)$ denotes the spectrum of $\Delta$.
The facts that $\sigma(\Delta)=[-4,0]$ and $\sigma(\Delta)$ is
purely continuous will be shown in Section 3.

The discrete Laplacian operator is one of the most important and
oldest difference operator [D] which is closely related to the
operator I'm studying in this paper. But before I introduce this
operator, I would like to talk about Schr\"{o}dinger's equation
and Schr\"{o}dinger operators.

A Schr\"{o}dinger's equation without a potential term is a partial
differential equation defined as $i\frac{\partial u}{\partial
t}=-\Delta u$ [E]. On the other hand,  Schr\"{o}dinger operators
acting on $L^2(\mathbb{R})$ are operators such as $Hf(x)=-\Delta
f(x)+v(x)f(x)$, where $V$ is a real-valued function on
$\mathbb{R}$ and is called a potential [D]. We may impose
condition on $V$ such as choosing $V \in L^{1}(\mathbb{R})$ so
that $H$ is an self-adjoint operator on $L^{2}(\mathbb{R})$. For a
detailed proof of this specific condition, see Davies [D].
Moreover, `$f \in L^2(\mathbb{R})$ with $\|f\|_{2}=1$ is called a
wave packet or state, and represents the instantaneous
configuration of a collection of electrons, atoms and molecules.
The operator $H$ is also called the Hamiltonian for historical
reasons - quantum theory can be regraded as a non-commutative
version of classical Hamiltonian mechanics. The evolution of a
quantum system is controlled by the Schr\"{o}dinger's equation
$i\frac{\partial f}{\partial t}=Hf$ with solution
$f(x,t)=e^{iHf}f(x,0)$' [D].


Finally, I want to now introduce the problem of my project. The
term `$\Delta$' which you will see below is a coefficient, not the
discrete Laplacian operator mentioned above.

For a particular example of a Schr\"{o}dinger's equation which is
also the main subject of my thesis, a time-dependent discrete
Schr\"{o}dinger's equation defined on $\ell ^{2}(\mathbb{Z})$ can
be written as follows:

\begin{equation*}
i\frac{d}{dt}\nu_{x}=-\frac{1}{\Delta}(\nu_{x-1}+\nu_{x+1})+\epsilon_{x}\nu_{x}
\end{equation*}
where
\begin{equation*}
\epsilon_{x}=\frac{2\cosh(\eta(x-r))^{2}}{\cosh(\eta(x-1-r))\cosh(\eta(x+1-r))}
\end{equation*}
with $\Delta=\cosh(\eta)$, $\eta\in \mathbb{R}^{+}$,
$x\in\mathbb{Z}$, $r\in\mathbb{R}$.

In order to understand the operator $i\frac{d}{dt}$, we define and
study $\mathcal{H}:\ell ^{2}(\mathbb{Z}) \rightarrow \ell
^{2}(\mathbb{Z})$ as follows

\begin{equation*}
\mathcal{H}=-\frac{1}{\Delta}(\mathcal{S}+\mathcal{S^{*}})+\epsilon_{x}\mathcal{I}
\end{equation*}

\noindent such that
\[
\ \mathcal{H}=\left(
\begin{array}{ccccccccc}
\ddots & \ddots & \ddots & \vdots &  & \vdots &  & \vdots &  \\
\ddots & \epsilon_{-x} & -1/\Delta & 0 & \ldots & 0 & \ldots & 0 & \ldots \\
\ddots & -1/\Delta & \ddots & \ddots & \ddots & \vdots &  & \vdots &  \\
\ldots & 0 & \ddots & \epsilon_{-1} & -1/\Delta & 0 & \ldots & 0 & \ldots \\
 & \vdots & \ddots & -1/\Delta & \epsilon_{0} & -1/\Delta & \ddots & \vdots &  \\
 \ldots & 0 & \ldots & 0 & -1/\Delta & \epsilon_{1} & \ddots & 0 & \ldots \\
 & \vdots &  & \vdots & \ddots & \ddots & \ddots & -1/\Delta & \ddots \\
 \ldots & 0 & \ldots & 0 & \ldots & 0 & -1/\Delta & \epsilon_{x} & \ddots \\
 & \vdots &  & \vdots &  & \vdots &  & \ddots & \ddots \\
\end{array}
\right)
\]

$\mathcal{H}$ is linear and bounded. The fact that $\mathcal{H}$
is linear is trivial. $\mathcal{H}$ is bounded because
$\|\mathcal{H}\|\leq
\|\frac{1}{\Delta}\mathcal{S}\|+\|\frac{1}{\Delta}\mathcal{S^{*}}\|+\|\epsilon{.}\mathcal{I}\|
\leq \|\mathcal{S}\|+\|\mathcal{S^{*}}\|+\|2\mathcal{I}\|=1+1+2=4$
since the maxima of $1/\Delta$ and $\epsilon_{x}$ are $1$ and $2$,
respectively. The reason that the maximum of $\epsilon_{x}$ is $2$
for any $r$ and $\eta$ is because of the properties of
$\epsilon_{x}$ which will explained in Section 2, and also Figure
1 below plots an example of $\epsilon_{x}$ versus $x$.

As you can see that the two main differences between the discrete
Laplacian operator and $\mathcal{H}$ are the tridiagonal terms of
their matrices where those terms in the discrete Laplacian
operator are constant functions: $1$, $-2$, and those in
$\mathcal{H}$ are variables depending on $r$ and $\eta$:
$-1/\Delta$, $\epsilon_{x}$, respectively.

We know that the spectrum of the discrete Laplacian operator
contains only the continuous spectrum when the diagonal terms of
the operator are constants $-2$. On the other side, the diagonal
terms $\epsilon_{x}$ of the operator $\mathcal{H}$ have a trapping
region around their centers and are $2$ for large $x$. Figure 1
presented below shows an example of $\epsilon_{x}$ with certain
parameters $r$ and $\eta$. So, $\mathcal{H}$ is a perturbed case
of the discrete Laplacian operator, and the nature of the
structure of the spectrum of $\mathcal{H}$ is different from it of
the discrete Laplacian operator because of the term
$\epsilon_{x}$.

\begin{figure}[h]
\includegraphics[width=10cm]{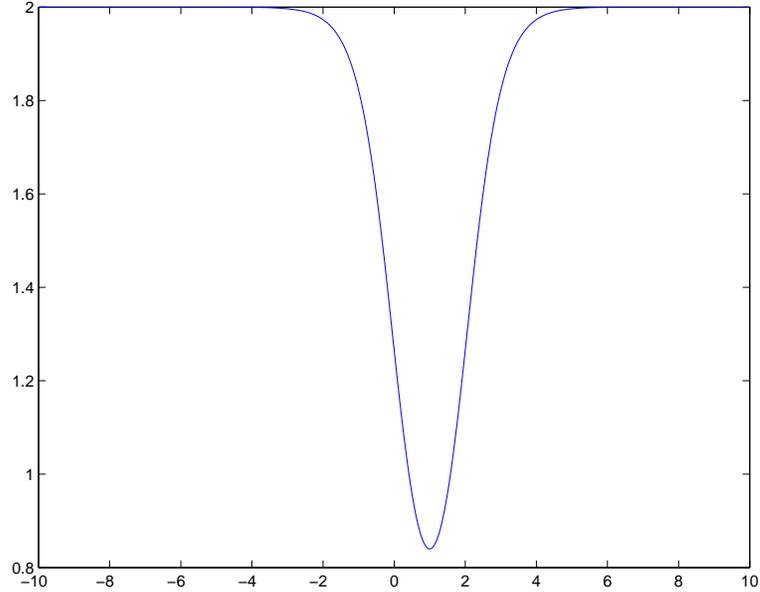}
\centering\caption{$\epsilon_{x}$ vs $x$ with $r=1$ and $\eta=1$}
\end{figure}

In this paper, I study the spectrum of $\mathcal{H}$ defined above
based on the discrete Laplacian operator which is an unperturbed
case of $\mathcal{H}$ and also some numerical results that can be
used to investigate the solutions of the time-dependent discrete
Schr\"{o}dinger's equation from the following three aspects:
\begin{itemize}
    \item stationary solution ($ \mathcal{H}\nu=0 $)
    \item periodic solution ($ \mathcal{H}\nu=\lambda\nu $)
    \item properties or representatives of the general solution
\end{itemize}

I have only studied the stationary and periodic solutions but not
yet the properties or representatives of general solution here.

For infinite dimensional $\mathcal{H}$ depending on $r$ and
$\eta$, my research shows that there exist three eigenvalues
depending on a certain region of $r$ and $\eta$. The first
eigenvalue of $\mathcal{H}$ which equals zero is exact. A detailed
analysis regrading this eigenvalue also called the zero mode [MN]
will be discussed in Section 3. Moreover, a second eigenvalue
exists which is proved by Michoel and Nachtergaele [MN]. Finally,
we also predict that a third eigenvalue exists for certain
parameters $r$ and $\eta$ based my numerical results. The primary
focus in my thesis is to find critical values or intervals of the
parameters $r$ and $\eta$ where these eigenvalues exist by
numerics. More precisely, the main object is to determine the
particular regions where there exist one, two, or three
eigenvalues for infinite dimensional $\mathcal{H}$. In order to
study this operator by numerics, we concentrate on a finite
dimensional subspace of $\ell^{2}(\mathbb{Z})$ such that
$\mathcal{H}$ is a finite dimensional operator. Since there always
exists a set of eigenvalues for finite dimensional $\mathcal{H}$,
we therefore focus on the eigenvalues that are isolated from the
set of all other eigenvalues which will belong to the continuous
spectrum and consider those isolated ones as the candidates of
eigenvalues which will belong to the point spectrum for infinite
dimensional $\mathcal{H}$. Those `isolated eigenvalues' I
mentioned can be seen more clearly on Figure 3 and Figure 4 in
Section 3.

My numerical results are created by three programs written by
myself in MATLAB. First, I graph the spectrum of $\mathcal{H}$
versus both parameters $r$ and $\eta$ separately to see how the
set of the point spectrum is distributed with different $r$ and
$\eta$. Second, I graph the eigenvectors with respect to the three
eigenvalues versus their components to see how the eigenvectors
are distributed and I'm hoping to see that the eigenvectors are
still nonzero for large matrix size $n$ with certain $r$ and
$\eta$. Third, I calculate the maxima of the square normalized
components of all those three eigenvectors. I compare the maxima
with the same parameters and different matrix sizes individually
and my expectation is that if the maxima do not decrease at
certain values or intervals of $r$ and $\eta$, then the
eigenvectors are not zero which also implies that there exist
eigenvalues for infinite dimensional $\mathcal{H}$. This step
mainly support my results. Finally, for the third eigenvalue
mentioned above, I also observe three components of its
corresponding eigenvector to see the decreasing rates of them. My
meaning of decreasing rates will be explained in the next few
sections. This process is to determine whether my conclusion of
the existence of the third eigenvalue at certain $r$ and $\eta$
based on the previous steps are correct or not. The next three
sections show my analysis in detail of the problem and the last
section concludes my entire paper, and they are presented in the
order as follows:
\begin{description}
    \item[Primary Results] This section states the procedures and the main
    discoveries in my project.
    \item[The Spectrum] This section discusses the definition of spectrum and also the structure and feature of the spectrum of
    $\mathcal{H}$.
    \item[The Eigenvalues] This section shows my numerical
    results for all the eigenvalues that are found.
    \item[Conclusion] This section summarizes my results and also
    talks about some further discussions of this operator.
\end{description}
\newpage

\section{Preliminary Results}

If $\mathcal{H}$ is a finite dimensional operator, then it is also
a compact operator. A compact operator is bounded, and a compact
operator that is also symmetric is self-adjoint. For more details,
see Hunter and Nachtergaele [HN]. Therefore, $\mathcal{H}$ is a
self-adjoint operator. Moreover, $\mathcal{H}$ is a nonnegative
operator. In order to prove this fact, we need to use the
following lemma first.

\begin{lemma}
Let $\mathcal{A}$ be a $2 \times 2$ matrix such that
\[
\ \mathcal{A}=\left(
\begin{array}{cc}
  a & c \\
  \overline{c} & b \\
\end{array}
\right)
\]
where $a,b\in\mathbb{R}$, $c\in\mathbb{C}$. Then $\mathcal{A} \geq
0$ if and only if $a,b \geq 0$ and $ab \geq |c|^{2}$.
\end{lemma}

\textbf{Proof.} Let $\mathcal{A} \geq 0$. First, $\mathcal{A} \geq
0$ if and only if all its eigenvalues $\lambda_{1},\lambda_{2}
\geq 0$. Second, the determinant of $\mathcal{\mathcal{A}}$ after
simplification is $\lambda^{2}-(a+b)\lambda+ab-|c|^{2}=0$.
Moreover,
\[ \left\{ \begin{array}{c}
\lambda_{1}+\lambda_{2}=a+b \\
\lambda_{1}\lambda_{2}=ab-|c|^2
\end{array}
\right. \]

Since $\lambda_{1},\lambda_{2} \geq 0$, we have $ab \geq |c|^{2}
\geq 0$ and $a,b \geq 0$. Conversely, if $a,b \geq 0$ and $ab \geq
|c|^{2}$, then $\lambda_{1},\lambda_{2} \geq 0$ which also implies
that $\mathcal{A} \geq 0$.

Now, we can decompose $\mathcal{H}$ into infinitely many
submatrices such that $\mathcal{H}=\sum_{x\in \mathbb{Z}}h_{x}$
where

\[
\ h_{x}=\left(
\begin{array}{cccccc}
  \ddots & \ddots & \vdots &  & \vdots &  \\
  \ddots & 0 & 0 & \ldots & 0 & \ldots \\
  \ldots & 0 & \epsilon_{x}^{+} & -1/\Delta & \vdots &  \\
   & \vdots & -1/\Delta & \epsilon_{x+1}^{-} & 0 & \ldots \\
  \ldots & 0 & \ldots & 0 & 0 & \ddots \\
   & \vdots &  & \vdots & \vdots & \ddots \\
\end{array}
\right)
\]

because $\epsilon_{x}=\epsilon_{x}^{+}+\epsilon_{x}^{-}$ where
\begin{equation*}
\epsilon_{x}^{\pm}=\frac{\cosh(\eta(x-r))}{\Delta\cosh(\eta(x\pm1-r))}
\end{equation*}

Since only four entries in $h_{x}$ are nonzero, if we want to
prove that $h_{x}$ is nonnegative, it is sufficient to focus on a
$2 \times 2$ matrix $\widetilde{h_{x}}$ defined as follows:

\[
\ \widetilde{h}_{x}=\left(
\begin{array}{cc}
  \epsilon_{x}^{+} & -1/\Delta \\
  -1/\Delta & \epsilon_{x+1}^{-} \\
\end{array}
\right)
\]

Since $\epsilon_{x}^{\pm} \geq 0$ for all $x\in \mathbb{Z}$,
$\epsilon_{x}^{+} \geq 0$ and $\epsilon_{x+1}^{-} \geq 0$. Thus,
the remaining is to show that whether
$\epsilon_{x}^{+}\epsilon_{x+1}^{-} \geq (1/\Delta)^{2}$ is true
or not. The calculation is presented below.

\begin{eqnarray*}
\epsilon_{x}^{+}\epsilon_{x+1}^{-} - (1/\Delta)^{2} & = & (\frac{\cosh(\eta(x-r))}{\Delta\cosh(\eta(x+1-r))}) \times \nonumber\\
 &  & (\frac{\cosh(\eta(x+1-r))}{\Delta\cosh(\eta(x-r))})-(1/\Delta)^{2}\\
 & = & 0
\end{eqnarray*}

Hence, $\widetilde{h}_{x}$ is nonnegative for all $x\in
\mathbb{Z}$ because $\epsilon_{x}^{+},\epsilon_{x+1}^{-} \geq 0$
and $\epsilon_{x}^{+}\epsilon_{x+1}^{-} = (1/\Delta)^{2}$. So,
$h_{x}$ is nonnegative for all $x\in \mathbb{Z}$ which implies
that $\mathcal{H}$ is nonnegative and also all eigenvalues of
$\mathcal{H}$ are nonnegative.

For infinite dimensional $\mathcal{H}$, my numerical results
predict that there exist three isolated eigenvalues of
$\mathcal{H}$ for a certain range of the parameters $r$ and $\eta$
such that the point spectrum of $\mathcal{H}$ is not empty, unlike
the discrete Laplacian operator. My conclusions are mainly based
on calculating the maxima of the square normalized components of
all the three eigenvectors corresponding to their eigenvalues with
various $r$ and $\eta$ and also observe the pictures of all the
three eigenvectors and the decreasing rates for three chosen
elements of the eigenvector of the third eigenvalue to compare
them with my numerics because the range of the existence of the
third isolated eigenvalue with respect to $r$ and $\eta$ is
different from it of the other two isolated eigenvalues. Roughly
speaking, with certain interval of $\eta$, the first two
eigenvalues exist for all $r$ in a certain region which will be
determined but the third eigenvalue does not exist for all $r$ in
that region. More details are presented in the following few
paragraphs.

I am using finite dimensional $\mathcal{H}$ to approximate
infinite dimensional $\mathcal{H}$. The reason why I follow such
path is based on what we know about the discrete Laplacian
operator. If $\eta=0$, then $\mathcal{H}$ is different from the
discrete Laplacian operator only by a sign. Moreover, my numerical
data includes the case where $\eta=0$ and it shows that the
eigenvectors approach to zero when the size of $\mathcal{H}$
becomes larger and larger though not presented here. So, we expect
that all the eigenvectors of $\mathcal{H}$ when $\eta=0$ converge
to zero by the $\sup$ norm as the size of $\mathcal{H}$ goes to
infinity which indicates that the limits of the eigenvalues are in
fact not eigenvalues. The reason that the converge is in the sense
of $\sup$ norm is because my numerical data shows that the maxima
of the square normalized components converge to zero, and so the
convergence is also in the weak $\ell^{2}$ sense. Moreover, we
already know that the spectrum of the discrete Laplacian operator
does not contain any eigenvalue. Therefore, using finite cases of
$\mathcal{H}$ to predict infinite cases of $\mathcal{H}$ is
reliable even if we change the parameters $r$ and $\eta$.

Now I want to talk about why I start with certain $r$ and $\eta$
out of other possibilities. For an n-dimensional $\mathcal{H}$
where $n\in\mathbb{N}$, the diagonal terms
$(\epsilon_{x})_{x=-n}^{n}$ are exponentially localized functions.
$\epsilon_x$ depends on $r$ and $\eta$ in which $r$ plays the role
of shifting $\epsilon_x$ and $\eta$ concentrates $\epsilon_x$
around its center, and $\epsilon_x$ is also symmetric around its
center shown in Figure 1 above. It is therefore sufficient to
focus on the case $r\in[n/2,n/2+1]$ where $n$ represents the size
of $\mathcal{H}$. Moreover, $\epsilon_x$ is an even function
because it is governed by hyperbolic cosine functions which are
even. As a result, we can begin with $\eta=0$ and then gradually
increase $\eta$. By plotting the graphs for the set of eigenvalues
of $\mathcal{H}$ versus the parameter $r\in[n/2,n/2+1]$ with
increasing values of $\eta$ which will be displayed in the next
section, the range of the set of the spectrum of $\mathcal{H}$
becomes smaller and smaller and three isolated eigenvalues are
extracted from the spectrum. For infinite dimensional case, these
three eigenvalues will belong to the point spectrum and all others
will belong to the continuous spectrum. The concept and structure
of spectrum will also be discussed in the following section.

From my numerical results, there always exist two eigenvalues for
all $r\in[n/2,n/2+1]$ and $\eta>0$. Moreover, for
$(r-[r])\in[0,0.1]\cup[0.9,1]$ and $\eta>1$, and for
$(r-[r])\in[0.2,0.3]\cup[0.7,0.8]$ and $\eta>2$, there exist a
total of three eigenvalues. The graph shown below predicts the
general results for the locations of the existence of the three
isolated eigenvalues with different $r$ and $\eta$ for infinite
approximation of $\mathcal{H}$. For numerical data which is not
completely presented here, the third eigenvalue exists when
$(r-[r]) \approx 0,1$ and $\eta=1$; $(r-[r]) \approx 0.15,0.85$
and $\eta=1.5$; $(r-[r]) \approx 0.25,0.75$ and $\eta=2$, and the
two curves in the graph are produced based on connecting those
points. More numerical results are therefore needed to verify this
graph.

\begin{figure}[h]
\includegraphics[width=10cm]{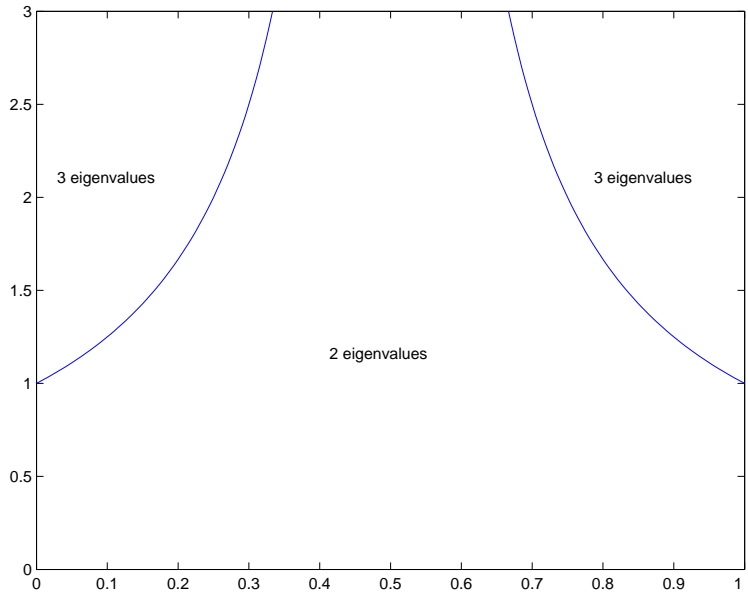}
\centering\caption{$\eta$ vs $(r-[r])$, where $[r]$ denotes the
integral part of $r$}
\end{figure}

The results above show that the point spectrum of infinite
dimensional $\mathcal{H}$ is nonempty at these certain regions.
Hunter and Nachtergaele [HN] and Michoel and Nachtergaele [MN]
have proved the existence of the first and second isolated
eigenvalues analytically, respectively. My results mainly come
from comparing the maxima of the square components of normalized
eigenvectors with different matrix size $n$ at certain values of
$r$ and $\eta$ numerically. The idea is that when the maxima do
not decrease for certain $r$ and $\eta$ while the matrix size $n$
becomes larger and larger, the eigenvalue does not vanish because
its corresponding eigenvector is nonzero. More details of my
numerical calculations of the results presented above will be
discussed in the next two sections.

My numerical results including the maxima and the graphs of
eigenvectors are created by a program written by myself in MATLAB
named "mainresults.m" and the one for observing the decreasing
rates is called "decreasingrates.m". Also, the pictures of the
spectrum of $\mathcal{H}$ versus $(r-[r])$ and $1/\Delta$,
respectively, are generated by "spectrum.m". The purposes of these
three programs will be explained in Section 3 and Section 4. The
next section explains the spectrum of $\mathcal{H}$ in more
details and also discusses the spectrum of the discrete Laplacian
operator.
\newpage

\section{The Spectrum}

In some circumstances, there are no eigenvalues for an infinite
dimensional bounded linear operator and so it is not possible to
expect to find an orthonormal basis that consists entirely of
eigenvectors. Thus, we need to define the spectrum in a more
general way, instead of considering it only contains eigenvalues
[HN].

The following two definitions are rewritten from `Applied
Analysis' by Hunter and Nachtergaele [HN].

\begin{definition}
Let $\mathcal{A}$ be a bounded operator defined on an infinite
dimensional Hilbert space. The resolvent set of $\mathcal{A}$,
denoted by $\rho(\mathcal{A})$, is the set of complex numbers
$\lambda$ such that $(\mathcal{A}-\lambda \mathcal{I})$ is
one-to-one and onto. The spectrum of $A$, denoted by
$\sigma(\mathcal{A})$, is the complement of the resolvent set in
$\mathcal{C}$, meaning that
$\sigma(\mathcal{A})=\mathbb{C}\backslash\rho(\mathcal{A})$.
\end{definition}

According to the open mapping theorem, $(\mathcal{A}-\lambda
\mathcal{I})^{-1}$ is bounded if $\mathcal{A}-\lambda \mathcal{I}$
is one-to-to and onto, and therefore $\mathcal{A}-\lambda
\mathcal{I}$ and $(\mathcal{A}-\lambda \mathcal{I})^{-1}$ are
one-to-one, onto, bounded linear operators when $\lambda \in \rho
(\mathcal{A})$. The following definition gives the structure of
the spectrum of a bounded linear operator.

\begin{definition}
Suppose that $\mathcal{A}$ is a bounded linear operator on a
Hilbert space.
\begin{enumerate}
    \item The point spectrum of $\mathcal{A}$ consists of all
    $\lambda\in\sigma(\mathcal{A})$ such that $\mathcal{A}-\lambda
    \mathcal{I}$ is not one-to-one. In this case $\lambda$ is
    called an eigenvalue of $\mathcal{A}$
    \item The continuous spectrum of $\mathcal{A}$ consists of all
    $\lambda \in \sigma(\mathcal{A})$ such that $\mathcal{A}-\lambda
    \mathcal{I}$ is one-to-one but not onto, and the range of $\mathcal{A}-\lambda
    \mathcal{I}$ is dense in this Hilbert space.
    \item The residual spectrum of $\mathcal{A}$ consists of all
    $\lambda \in \sigma(\mathcal{A})$ such that $\mathcal{A}-\lambda
    \mathcal{I}$ is one-to-one but not onto, and the range of $\mathcal{A}-\lambda
    \mathcal{I}$ is not dense in this Hilbert space.
\end{enumerate}
\end{definition}

We define $\mathcal{H}$ on a finite subspace of
$\ell^{2}(\mathbb{Z})$ for our problem and we then have
$\mathcal{H}\nu=\lambda \nu$ where $\lambda$ is an eigenvalue and
$\nu$ is a nonzero eigenvector. Thus, the kernel of $\mathcal{H}$
does not contain zero, and so $\mathcal{H}-\lambda \mathcal{I}$ is
not one-to-one and $\lambda \in \sigma(\mathcal{H})$.

Because of the diagonal term $\epsilon_{x}$ in $\mathcal{H}$, the
structure of the spectrum of $ \mathcal{H} $ is different from it
of the discrete Laplacian operator. As shown in Figure 1 above,
the function $\epsilon_{x}$ remains constant for large $|x|$ but
drops toward zero for small $|x|$ around its center. The existence
of this dropping region makes the spectrum of $\mathcal{H}$ so
special in which the point spectrum is not empty such that we
expect there exist eigenvalues when $\mathcal{H}$ is infinite
dimensional. The location, width, and depth of this trapping
region depend significantly on the two parameters $r$ and $\eta$.

Since $\epsilon_{x}$ is even and symmetric, $r$ shifts
$\epsilon_{x}$ to the left or to the right depending on its sign,
and $\eta$ concentrates $\epsilon_{x}$ narrow or wide around the
center of $\epsilon_{x}$ depending on its magnitude, the most
interesting case for $\mathcal{H}$ is when the difference between
$x$ and $r$ in $\epsilon_{x}$ is small. It is therefore sufficient
to consider $x \in \mathbb{Z^{+}}$ and $r \in [n/2,n/2+1]$ where
$n$ is the dimension of $\mathcal{H}$. Moreover, we can choose
$\eta\geq 0$ because $\epsilon_{x}$ is an even function.

The following five pictures plot the point spectrum or eigenvalues
denoted them as $\lambda$ of a finite $\mathcal{H}$ versus
$r\in[n/2,n/2+1]$ with five different $\eta$ follows:

\begin{figure}[h]
\includegraphics[width=10cm]{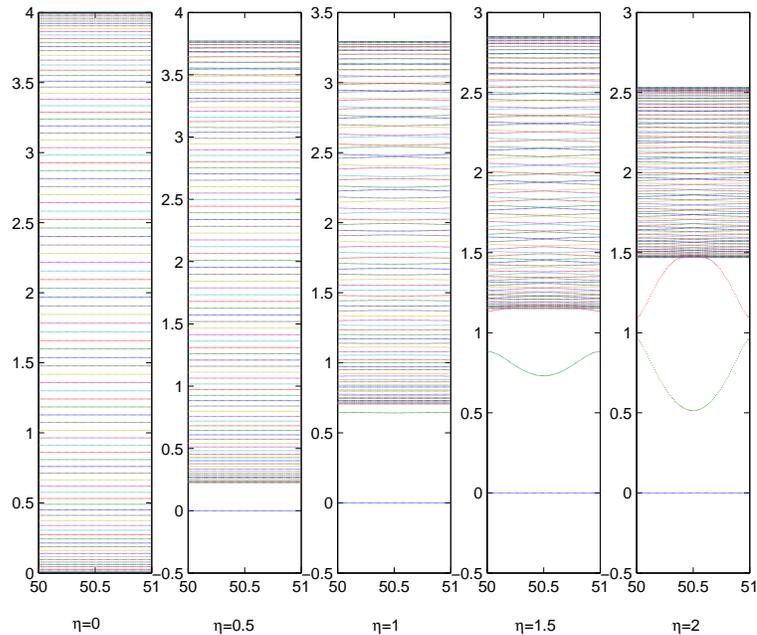}
\centering\caption{$\lambda$ vs $r$ with $n=100$ and
$r\in[n/2,n/2+1]$}
\end{figure}

The first graph $\eta=0$ shows that the range of the spectrum of
$\mathcal{H}$ is in $[0,4]$ which is similar as the discrete
Laplacian operator where it is in $[-4,0]$. The second graph
$\eta=0.5$ shows that there is one isolated eigenvalue $\lambda=0$
and the range of the other eigenvalues of $\mathcal{H}$ is smaller
than it when $\eta=0$. The third graph where $\eta=1$ and the
fourth graph where $\eta=1.5$ show that there are two isolated
eigenvalues and the range of all the other eigenvalues becomes
smaller and smaller even than before. At last, the fifth graph
shows that there are three isolated eigenvalues and the range of
the eigenvalues besides those three eigenvalues of $\mathcal{H}$
is about [1.4,2.6]. However, at $r=n/2+0.5$ the third appearing
eigenvalue seems to be still connecting with the spectrum, not
totally isolated. Before I draw initial conclusions of critical
values or intervals for $r$ and $\eta$, I like to show another
five pictures that plot the eigenvalues $\lambda$ of $\mathcal{H}$
versus $1/\cosh \eta$ where $\eta \in [0,2]$ with five different
$r$ as follows:

\begin{figure}[h]
\includegraphics[width=10cm]{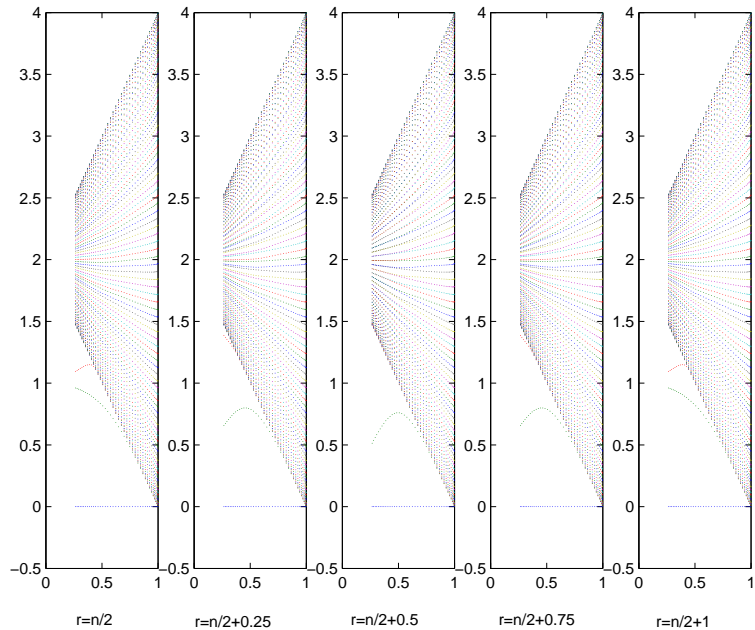}
\centering\caption{$\lambda$ vs $1/\cosh \eta$ with $ n=100 $ and
$\eta\in[0,2]$}
\end{figure}

The first graph where $r=n/2$ and the fifth graph where $r=n/2+1$
show that there are three isolated eigenvalues at $1/\cosh \eta
\approx 0.3$ and two isolated eigenvalues at $1/\cosh \eta \approx
0.6$. $0.3$ corresponds approximately to $\eta=2$ and $0.6$ is
about $\eta=1$. Moreover, the second graph where $r=n/2+0.25$ and
the fourth graph where $r=n/2+0.75$ show that there are still
three isolated eigenvalues when $1/\cosh \eta \approx 0.3$ and two
isolated eigenvalues when $1/\cosh \eta \approx 0.6$, but now the
third eigenvalue is not quite isolated from the set of other
eigenvalues and seems to merge into the other sets of point
spectrum. The third graph where $r=n/2+0.5$ show that there only
exist two isolated eigenvalues. So, by Figure 3 and Figure 4
above, we can see that the third appearing eigenvalue does not
seem to exist for all $r \in [n/2,n/2+1]$.

Therefore, we expect that the first and second isolated
eigenvalues appear somewhere at $\eta\in (0,1)$ and exist for all
$r \in [n/2,n/2+1]$, and the third isolated eigenvalue appears
somewhere at $\eta\in (1,2)$ and exist for for all $r \in
[n/2,n/2+0.5)\cup (n/2+0.5,n/2+1]$. By approximating infinite
$\mathcal{H}$ from finite $\mathcal{H}$, these three eigenvalues
will remain in the point spectrum and all other eigenvalues will
belong in the continuous spectrum.

You may ask that how numerical calculations for finite size
systems can give insight in the infinite size systems. Since
$\mathcal{H}$ is defined on $\ell^{2}(\mathbb{Z})$, we can focus
on a subspace of $\ell^{2}(\mathbb{Z})$ such that $\mathcal{H}$ is
finite dimensional. Eigenvalues and eigenvectors of finite
$\mathcal{H}$ can be calculated by numerics and we can increase
the size of $\mathcal{H}$ since the eigenvectors are defined on
$\ell^{2}(\mathbb{Z})$. The spectrum of infinite dimensional
$\mathcal{H}$ can then be predicted by observing it of finite
dimensional $\mathcal{H}$.

Moreover, we know the structure of the spectrum of the discrete
Laplacian operator and $\mathcal{H}$ is a perturbed case of it, we
can study the spectrum of the discrete Laplacian operator in order
to understand the spectrum of $\mathcal{H}$ based on a theorem
called the compact perturbation theorem below.

\begin{theorem}
Let $\mathcal{T}\in B(\mathbf{X})$ and let $\mathcal{A}$ be
$\mathcal{T}-compact$. Then $\mathcal{T}$ and
$\mathcal{T}+\mathcal{A}$ have the same essential spectrum.
\end{theorem}

This theorem is rewritten from `Perturbation Theory for Linear
Operators' by Kato [K].

$\mathcal{H}$ is constructed basically by reversing the sign of
the discrete Laplacian operator and adding a compact perturbation
term $\epsilon_x$. The operator $\epsilon_{x}\mathcal{I}$ in
$\mathcal{H}$ is in fact a compact perturbation. To prove this
statement, we can first define
$\mathcal{T}=(\epsilon_{x}-2)\mathcal{I}$ because we intend to
show that $\epsilon_{x}\mathcal{I}$ is perturbed from
$2\mathcal{I}$ which is the diagonal of the discrete Laplacian
operator with inverse sign. We can then define an operator
$\mathcal{T}_{n}$ with finite-dimensional range such that
$\mathcal{T}_{n}$ is compact [HN] and

\[
\ \mathcal{T}_{n}=\left(
\begin{array}{ccccccccccc}
  \ddots & \vdots & \vdots &  &  & \vdots & \vdots &  &  & \vdots &  \\
  \ddots & 0 & 0 & \ldots & \ldots & 0 & 0 & \ldots & \ldots & 0 & \ldots \\
  \ldots & 0 & \epsilon_{-x} & \ddots &  & \vdots & \vdots &  &  & \vdots &  \\
   & \vdots & \ddots & \ddots & \ddots & \vdots & \vdots &  &  & \vdots &  \\
   & \vdots &  & \ddots & \epsilon_{-1} & 0 & 0 & \ldots & \ldots & 0 & \ldots \\
  \ldots & 0 & \ldots & \ldots & 0 & \epsilon_{0} & 0 & \ldots & \ldots & 0 & \ldots \\
  \ldots & 0 & \ldots & \ldots & 0 & 0 & \epsilon_{1} & \ddots &  & \vdots &  \\
   & \vdots &  &  & \vdots & \vdots & \ddots & \ddots & \ddots & \vdots &  \\
   & \vdots &  &  & \vdots & \vdots &  & \ddots & \epsilon_{x} & 0 & \ldots \\
  \ldots & 0 & \ldots & \ldots & 0 & 0 & \ldots & \ldots & 0 & 0 & \ddots \\
   & \vdots &  &  & \vdots & \vdots &  &  & \vdots & \vdots & \ddots \\
\end{array}
\right)
\]

If $\lim_{n\rightarrow\infty}\|\mathcal{T}_{n}-\mathcal{T}\|=0$,
then $\mathcal{T}$ is compact where $\|\cdot\|$ denotes the
operator norm. Now, we let $\nu\in\ell^{2}(\mathbb{Z})$ and

\begin{eqnarray*}
  \|(\mathcal{T}_{n}-\mathcal{T})\nu\|_{\ell^{2}(\mathbb{Z})} & \leq & (\sum_{|x|>n}(\epsilon_{x}-2)^{2}\nu_{x}^{2})^{1/2} \\
   & \leq & \max_{|x|>n}|\epsilon_{x}-2| (\sum_{|x|>n}\nu_{x}^{2})^{1/2}\\
   & = & \max_{|x|>n}|\epsilon_{x}-2| \|\nu\|_{\ell^{2}(\mathbb{Z})}
\end{eqnarray*}

So, $\lim_{n\rightarrow\infty}\|\mathcal{T}_{n}-\mathcal{T}\|=0$
because $\lim_{x\rightarrow\infty}|\epsilon_{x}-2|=0$ is as shown
in Figure 1.

Thus, $\epsilon_{x}\mathcal{I}$ in $\mathcal{H}$ is a compact
perturbation and $\mathcal{H}$ and the discrete Laplacian operator
have the same essential spectrum. The spectrum of the discrete
Laplacian operator is purely continuous and is contained in the
interval $[-4,0]$ and the point spectrum of finite $\mathcal{H}$
at $\eta=0$ is contained in $[0,4]$ based on my numerics. When
$\eta=0$, $\mathcal{H}$ is different from the discrete Laplacian
operator by a sign. Moreover, my numerical results presented in
the next section shows that at $\eta=0$ the maxima of the square
normalized components of all the eigenvectors of the three
isolated eigenvalues decrease as the matrix size $n$ of
$\mathcal{H}$ increase and we expect that the eigenvectors
converge to zero for infinite $n$. So, there are no eigenvalues
when $\eta=0$ for infinite dimensional $\mathcal{H}$ and they all
converge to the set of the continuous spectrum. This fact agrees
with the case of the discrete Laplacian operator. Therefore, when
we change $r$ and $\eta$ of a finite n-dimensional $\mathcal{H}$,
we expect that the three isolated eigenvalues will converge to the
set of the point spectrum and all other eigenvalues will converge
to the set of point spectrum as $n$ goes to infinity. Finally, the
residual spectrum of the discrete Laplacian operator and
$\mathcal{H}$ are empty because they are both bounded and
self-adjoint [HN].

The mathematical argument that determines the continuous spectrum
in the infinite chain limit of the discrete Laplacian operator is
presented below.

Let $\mathcal{F}: L^{2}(\mathbb{T}) \rightarrow
\ell^{2}(\mathbb{R})$ denote the Fourier transform and
$\mathcal{F^{*}}$ be the inverse Fourier transform. $\mathcal{F}$
is an unitary operator. For proofs, see Hunter and Nachtergaele
[HN]. Moreover, consider only now that
$\Delta=\mathcal{S}+\mathcal{S^{*}}-2\mathcal{I}$ is the discrete
Laplacian operator.

It follows that $\mathcal{F^{*}}\Delta\mathcal{F}=\mathcal{D}$
where $\mathcal{D}$ is a diagonal matrix that consists of spectral
elements of $\Delta$ on the diagonal since $\Delta$ is
self-adjoint. Moreover, if $f(x)\in L^{2}(\mathbb{T})$ where
$\mathbb{T}=[0,2\pi]$ and $\widehat{f_{n}}\in
\ell^{2}(\mathbb{Z})$, then

\begin{eqnarray*}
\mathcal{F^{*}}\Delta\mathcal{F}f(x) & = & \mathcal{F^{*}}(\mathcal{S}+\mathcal{S^{*}}-2\mathcal{I})\mathcal{F}f(x) \\
 & = & \mathcal{F^{*}}(\mathcal{S}+\mathcal{S^{*}}-2\mathcal{I})\widehat{f_{n}} \\
 & = & \mathcal{F^{*}}(\mathcal{S}\widehat{f_{n}}+\mathcal{S^{*}}\widehat{f_{n}}-2\mathcal{I}\widehat{f_{n}}) \\
 & = & \mathcal{F^{*}}(\widehat{f_{n-1}}+\widehat{f_{n+1}}-2\widehat{f_{n}}) \\
 & = & \frac{1}{\sqrt{2\pi}}\sum_{n\in \mathbb{Z}}(\widehat{f_{n-1}}+\widehat{f_{n+1}}-2\widehat{f_{n}})e^{inx} \\
 & = & \frac{1}{\sqrt{2\pi}}\sum_{n\in
 \mathbb{Z}}\widehat{f_{n-1}}e^{inx}+\frac{1}{\sqrt{2\pi}}\sum_{n\in
 \mathbb{Z}}\widehat{f_{n+1}}e^{inx}-\frac{1}{\sqrt{2\pi}}\sum_{n\in
 \mathbb{Z}}2\widehat{f_{n}}e^{inx} \\
 & = & (e^{ix})\frac{1}{\sqrt{2\pi}}\sum_{n\in \mathbb{Z}}\widehat{f_{n-1}}e^{i(n-1)x} +
(e^{-ix})\frac{1}{\sqrt{2\pi}}\sum_{n\in
 \mathbb{Z}}\widehat{f_{n+1}}e^{i(n+1)x} - \nonumber \\
 &   &(2)\frac{1}{\sqrt{2\pi}}\sum_{n\in \mathbb{Z}}\widehat{f_{n}}e^{inx} \\
 & = & (e^{ix}+e^{-ix}-2)f(x) \\
 & = & (2\cos(x)-2)f(x) \\
 & = & d(x)f(x) \\
 & = & \mathcal{D}f(x)
\end{eqnarray*}

It follows that
$\sigma(\Delta)=\sigma(\mathcal{D})=\{d(x)|x\in[0,2\pi]\}=[-4,0]$
and $\sigma(\Delta)$ is entirely continuous.

Moreover, the bottom of the continuous spectrum of $\mathcal{H}$
is given by a formula $2(1-\Delta^{-1})$ [MN] which agrees with my
numerics shown in Figure 3 above. The smallest eigenvalue, or the
first isolated eigenvalue of $\mathcal{H}$ in our case, is also
called the "ground state energy" and its corresponding
eigenfunction is called the "ground state" which is the
configuration of the system with the smallest total energy, and
the other eigenvalues correspond to discrete excitations of the
system [D].

As proved by Michoel and Nachtergaele [MN], the eigenfunction of
the first isolated eigenvalue of $\mathcal{H}$ is
$V_{x}^{(1)}=\frac{1}{\cosh (\eta(x-r))}\in \ell^{2}(\mathbb{Z})$.

The fact that this eigenfunction belongs to $\ell^{2}(\mathbb{Z})$
can proved as follows:
\begin{equation*}
\sum_{x\in \mathbb{Z}}|\frac{1}{\cosh(\eta(x-r))}|^{2}=\sum_{x\in
\mathbb{Z}}\frac{4}{{e^{2\eta(x-r)}+e^{-2\eta(x-r)}+2}}
\end{equation*}
and
\begin{equation*}
\frac{4}{{e^{2\eta(x-r)}+e^{-2\eta(x-r)}+2}} \leq
\frac{4}{e^{2\eta(x-r)}} = \frac{4e^{2\eta r}}{e^{2\eta x}}
\end{equation*}

Since $\sum_{x\in \mathbb{Z}}\frac{4e^{2\eta r}}{e^{2\eta x}}$ is
a geometric series, by comparison test, we can conclude that

\begin{equation*}
\sum_{x\in \mathbb{Z}}|\frac{1}{\cosh(\eta(x-r))}|^{2} < \infty
\end{equation*}

and so $V_{x}^{(1)}\in \ell^{2}(\mathbb{Z})$.

Moreover, the fact that $V_{x}^{(1)}$ is indeed an eigenvector of
$\mathcal{H}V_{x}^{(1)}=0$ for the zero mode can be verified as
follows:

\begin{eqnarray*}
\mathcal{H}(\frac{1}{\cosh(\eta(x-r))}) & = &
\epsilon_{x}(\frac{1}{\cosh(\eta(x-r))}) - \nonumber \\
 &   & \frac{1}{\Delta}(\frac{1}{\cosh(\eta(x-1-r))}+\frac{1}{\cosh(\eta(x+1-r))}) \\
 & = & (\frac{2\cosh(\eta(x-r))^{2}}{\cosh(\eta(x-1-r))\cosh(\eta(x+1-r))})(\frac{1}{\cosh(\eta(x-r))}) - \nonumber \\
 &   & \frac{1}{\cosh(\eta)}(\frac{1}{\cosh(\eta(x-1-r))}+\frac{1}{\cosh(\eta(x+1-r))}) \\
 & = & \frac{2\cosh(\eta(x-r))}{\cosh(\eta(x-1-r))\cosh(\eta(x+1-r))} - \nonumber \\
 &   & \frac{1}{\cosh(\eta)}(\frac{\cosh(\eta(x+1-r))+\cosh(\eta(x-1-r))}{\cosh(\eta(x-1-r))\cosh(\eta(x+1-r))}) \\
 & = & \frac{2\cosh(\eta(x-r))}{\cosh(\eta(x-1-r))\cosh(\eta(x+1-r))} - \nonumber \\
 &   & \frac{1}{\cosh(\eta)}(\frac{((2\cosh(\eta))(2\cosh(x-r)))/2}{\cosh(\eta(x-1-r))\cosh(\eta(x+1-r))}) \\
 & = & \frac{2\cosh(\eta(x-r))}{\cosh(\eta(x-1-r))\cosh(\eta(x+1-r))} - \nonumber \\
 &   & \frac{2\cosh(\eta(x-r))}{\cosh(\eta(x-1-r))\cosh(\eta(x+1-r))}\\
 & = & 0
\end{eqnarray*}

The following section explains my MATLAB programs, and it also
shows my numerical results for all the three isolated eigenvalues
in order.

\newpage

\section{The Eigenvalue Problem}

There exist two isolated eigenvalues of $\mathcal{H}$ when $\eta
> 0$ and $r\in [n/2,n/2+1]$. The third isolated eigenvalue appears
when $\eta > 1$ and the interval of $r$ where it does not exist
becomes smaller and smaller when $\eta$ becomes larger and larger.
Moreover, the third eigenvalue does not exist at $r=n/2+0.5$.
These numerical results are based on three programs which are
explained in the next paragraph.

The first program named "spectrum.m" plots the entire set of
eigenvalues of an n-dimensional $\mathcal{H}$ versus $r \in
[n/2,n/2+1]$ and $\frac{1}{\cosh \eta}$ where $\eta \in [0,2]$
separately as shown in Figure 3 and Figure 4 above. This program
gives the picture of the locations in which there exist isolated
eigenvalues. With various $r$ and $\eta$ in $\mathcal{H}$, the
second program named "mainresults.m" graphs the eigenvectors of
all the three isolated eigenvalues versus the components of the
eigenvectors, and this program also calculates the maxima of the
square normalized components in the eigenvectors of those
eigenvalues. The first part of "mainrestuls.m" shows the
distributions of the eigenvectors at different $r$ and $\eta$. The
second part compares the maxima described above to see when they
do not decrease for certain $r$ and $\eta$ as $n$ increases which
imply the eigenvectors are nonzero and the existence of
eigenvalues. The last program named "decreasingrates.m" is
designed only for the third appearing eigenvalue because it does
not exist for all $r\in [n/2,n/2+1]$. I claim that the third
eigenvalue does not exist when $r=n/2+0.5$ and so I investigate
the decreasing rate of three selected components of the
eigenvector of the third eigenvalue. We let $ \nu_{x}^{(n)}$
denote an n-dimensional eigenvector at components $x=n/4,n/2,3n/4$
and assume that $|\nu_{x}^{(n)}|^{2} \leq e^{-ax}$ because the
square normalized component is decreasing as shown on Figure 7
below. $a=-\frac{log(|\nu_{x}^{(n)}|^{2})}{x}$ is then defined to
be the decreasing rate. We are expecting to see that $a\rightarrow
0$ as $r \rightarrow (n/2+0.5)$.

The codes of all the programs are described in Appendix B below.
The following three subsections show some selected numerical
results for the three eigenvalues denoted them as $\lambda_{1}$,
$\lambda_{2}$, and $\lambda_{3}$, respectively.

\subsection{The first isolated eigenvalue}

The pictures below graph the eigenvector denoted as $V_{x}^{(1)}$
with respect to $\lambda_{1}$ versus the components denoted the
index as $x$ of $V_{x}^{(1)}$ when $n=100$. From top to bottom,
$r$ increases from $n/2$ to $n/2+1$ by $0.25$ each time, and from
left to right, $\eta$ increases from $0$ to $2$ by $0.5$ each
time.
\begin{figure}[h]
\includegraphics[width=10cm]{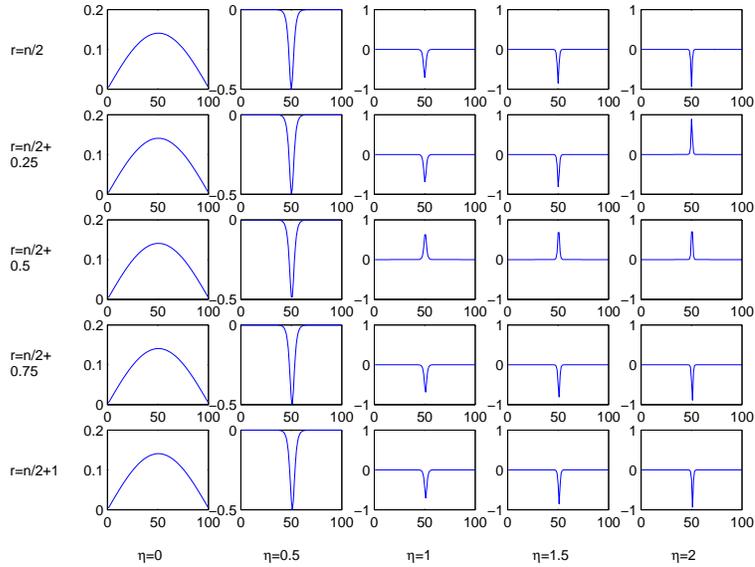}
\centering\caption{$V_{x}^{(1)}$ vs $x$ with $n=100$}
\end{figure}

As you can see from Figure 5 above, when $\eta > 0$, $V_{x}^{(1)}$
starts to concentrate around $x=n/2$ for all $r\in [n/2,n/2+1]$
and there exists an absolute maximum at $x=n/2$ on all the graphs.
This maximum shows that $V_{x}^{(1)}$ is nonzero although most
components of $V_{x}^{(1)}$ are zero and implies the existence of
the eigenvalue $\lambda_{1}$ for infinite $\mathcal{H}$. Similarly
for the cases of larger $n$ which are not shown here,
$V_{x}^{(1)}$ has an absolute maximum at $x=n/2$ and is zero
almost everywhere else.

Next, I want to justify that if the absolute maximum mentioned
above does not decrease for all $\eta > 0$ and $r \in [n/2,n/2+1]$
as $n$ increases.

The table below shows the maxima of the square normalized
components in $V_{x}^{(1)}$ at $n=200$ and $\eta \in [0.2]$ where
$\eta$ is increased by $(2-0)/4=0.5$ every time, and $r$ ranges
from $n/2$ to $n/2+1$ and is increased by $0.25$ each time.

\begin{center}
\begin{tabular}{|c|c|c|c|c|c|}
  \hline
   & $\eta=0$ & $\eta=0.5$ & $\eta=1$ & $\eta=1.5$ & $\eta=2$ \\
  \hline
  $r=n/2$      & 0.009950 & 0.250000 & 0.498981 & 0.723493 & 0.874100 \\
  $r=n/2+0.25$ & 0.009950 & 0.246134 & 0.470008 & 0.653751 & 0.788057 \\
  $r=n/2+0.5$  & 0.009950 & 0.235004 & 0.394028 & 0.464358 & 0.488313 \\
  $r=n/2+0.75$ & 0.009950 & 0.246134 & 0.470008 & 0.653751 & 0.788057 \\
  $r=n/2+1$    & 0.009950 & 0.250000 & 0.498981 & 0.723493 & 0.874100 \\
  \hline
\end{tabular}
\end{center}

\vspace{0.25 in} The table below shows these maxima at $n=300$ and
the rest is the same as above. \vspace{0.25 in}

\begin{center}
\begin{tabular}{|c|c|c|c|c|c|}
  \hline
   & $\eta=0$ & $\eta=0.5$ & $\eta=1$ & $\eta=1.5$ & $\eta=2$ \\
  \hline
  $r=n/2$      & 0.006644 & 0.250000 & 0.498981 & 0.723493 & 0.874100 \\
  $r=n/2+0.25$ & 0.006644 & 0.246134 & 0.470008 & 0.653751 & 0.788057 \\
  $r=n/2+0.5$  & 0.006644 & 0.235004 & 0.394028 & 0.464358 & 0.488313 \\
  $r=n/2+0.75$ & 0.006644 & 0.246134 & 0.470008 & 0.653751 & 0.788057 \\
  $r=n/2+1$    & 0.006644 & 0.250000 & 0.498981 & 0.723493 & 0.874100 \\
  \hline
\end{tabular}
\end{center}

\vspace{0.25 in} By comparing the maxima of the square normalized
components in $V_{x}^{(1)}$ on these two tables at the same values
of $r$ and $\eta$ with two different $n$, it shows that the maxima
remain unchanged besides when $\eta=0$. When $\eta=0$, the maxima
decrease as $n$ increases for all five different $r$. So the
maxima begin to be nondecreasing somewhere at $\eta \in (0,0.5)$.
My next step is to compare the maxima between $n=300$ and $n=400$
with $\eta \in [0,0.5]$ and increase $\eta$ by $(0.5-0)/4=0.125$
each time and keep the same $r$ as above. The following table
shows the procedure of how I narrow down the interval of $\eta$
with different $n_{1},n_{2}$ where $n_{1} < n_{2}$ and
$\eta_{1},\eta_{2}$ where $\eta \in [\eta_{1},\eta_{2}]$.
\vspace{0.25 in}
\begin{center}
\begin{tabular}{|c|c|c|c|}
  \hline
  $n_{1}$ & $n_{2}$ & $\eta_{1}$ & $\eta_{2}$ \\
  \hline
  200 & 300 & 0 & 2 \\
  300 & 400 & 0 & 0.5 \\
  400 & 500 & 0 & 0.125 \\
  500 & 600 & 0 & 0.0625 \\
  600 & 700 & 0 & 0.03125 \\
  \hline
\end{tabular}
\end{center}

\vspace {0.25 in} For the tables from $n=700$ to $n=1000$ which
are not shown here when $\eta \in [0,0.03125]$, the interval of
$\eta$ starting at zero becomes smaller and smaller beginning with
$r=n/2$ and $r=n/2+1$ where the maxima of the square normalized
components in $V_{x}^{(1)}$ do not decrease, followed by
$r=n/2+0.5$, and then $r=n/2+0.25$ and $r=n/2+0.75$. Since the
spectrum for infinite $\mathcal{H}$ is purely continuous when
$\eta=0$ because it is only different from the discrete Laplacian
operator by a sign and the point spectrum becomes nonempty for an
$\eta$ in $(0,0.3125)$ which is an open set, the critical value
for $\eta$ where $\lambda_{1}$ exists for infinite $\mathcal{H}$
is believed to be at $\eta=0$.
\subsection{The Second isolated eigenvalue}
$\lambda_{2}$ exists in a very similar fashion as $\lambda_{1}$.
The pictures below draw the eigenvector denoted as $V_{x}^{(2)}$
with respect to $\lambda_{2}$ versus the components denoted the
index as $x$ of $V_{x}^{(2)}$ when $n=100$. From top to bottom,
$r$ is increased from $n/2$ to $n/2+1$ by $0.25$ each time, and
from left to right, $\eta$ is increased from $0$ to $2$ by $0.5$
each time.

\begin{figure}[h]
\includegraphics[width=10cm]{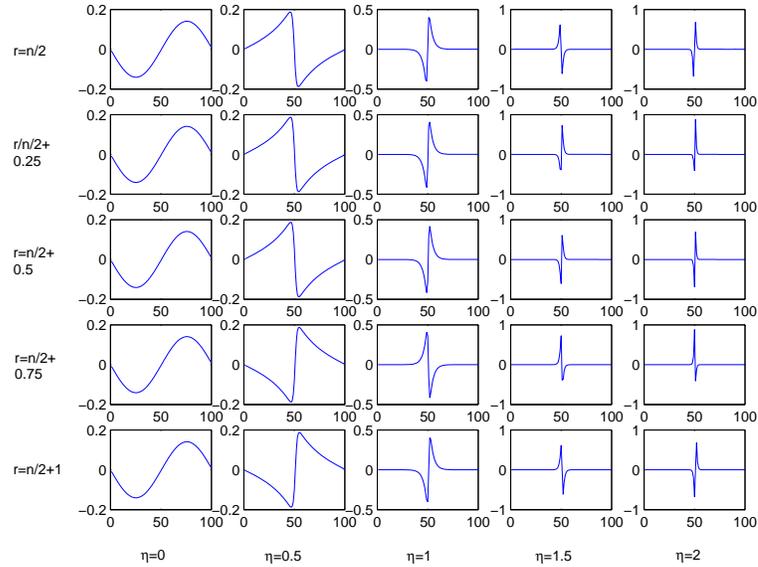}
\centering\caption{$V_{x}^{(2)}$ vs $x$ with $n=100$}
\end{figure}

As appeared on Figure 6 below, $V_{x}^{(2)}$ starts to be
concentrated obviously around $x=n/2$ like $V_{x}^{(1)}$ when
$\eta > 0.5$ for all $r\in [n/2,n/2+1]$ and there also exists an
absolute maximum at $x=n/2$ on those graphs. By graphing the same
pictures with larger $n$, there always exists an absolute maximum
in $V_{x}^{(2)}$. Thus, the eigenvector will remain nonzero and
$\lambda_{2}$ will not disappear for large $n$.

The two tables presented below show the maxima of the square
normalized components in $V_{x}^{(2)}$ at $n=200$ and $n=300$, and
$\eta \in [0.2]$ where $\eta$ increases by $(2-0)/4=0.5$ every
time, and $r$ ranges from $n/2$ to $n/2+1$ and increases by $0.25$
each time.

\begin{center}
\begin{tabular}{|c|c|c|c|c|c|}
  \hline
   & $\eta=0$ & $\eta=0.5$ & $\eta=1$ & $\eta=1.5$ & $\eta=2$ \\
  \hline
  $r=n/2$ & 0.009950 & 0.032095 & 0.159968 & 0.375804 & 0.459508 \\
  $r=n/2+0.25 $ & 0.009950 & 0.032057 & 0.171671 & 0.526199 & 0.772677 \\
  $r=n/2+0.5 $ & 0.009950 & 0.031823 & 0.172972 & 0.368510 & 0.476093 \\
  $r=n/2+0.75$ & 0.009950 & 0.032057 & 0.171671 & 0.526199 & 0.772677 \\
  $r=n/2+1$ & 0.009950 & 0.032095 & 0.159968 & 0.375804 & 0.459508 \\
  \hline
\end{tabular}
\end{center}

\vspace{0.25 in}

\begin{center}
\begin{tabular}{|c|c|c|c|c|c|}
  \hline
  & $\eta=0$ & $\eta=0.5$ & $\eta=1$ & $\eta=1.5$ & $\eta=2$ \\
  \hline
  $r=n/2$ & 0.006644 & 0.031971 & 0.159968 & 0.375804 & 0.459508 \\
  $r=n/2+0.25$ & 0.006644 & 0.031933 & 0.171671 & 0.526199 & 0.772677 \\
  $r=n/2+0.5$ & 0.006644 & 0.031699 & 0.172972 & 0.368510 & 0.476093 \\
  $r=n/2+0.75$ & 0.006644 & 0.031933 & 0.171671 & 0.526199 & 0.772677 \\
  $r=n/2+1$ & 0.006644 & 0.031971 & 0.159968 & 0.375804 & 0.459508 \\
  \hline
\end{tabular}
\end{center}

\vspace{0.25 in} The maxima shown on the tables above appeare to
be nondecreasing when $\eta > 0.5$ and my next step is to focus on
$\eta \in [0,1]$ and compare the maxima of the square normalized
components of $V_{x}^{(2)}$ between $n=300$ and $n=400$. My
complete process is shown in the table below: \vspace {0.25 in}

\begin{center}
\begin{tabular}{|c|c|c|c|}
  \hline
  $n_{1}$ & $n_{2}$ & $\eta_{1}$ & $\eta_{2}$ \\
  \hline
  200 & 300 & 0 & 2 \\
  300 & 400 & 0 & 1 \\
  400 & 500 & 0 & 0.75 \\
  500 & 600 & 0 & 0.5625 \\
  \hline
\end{tabular}
\end{center}

\vspace{0.25 in} When $n=600,700,800,900,1000$ for $\eta \in
[0,0.5625]$, same as in the previous situation but much slower,
the interval of $\eta$ starting at zero where the maxima do not
decrease becomes smaller and smaller beginning with $r=n/2$ and
$r=n/2+1$, followed by $r=n/2+0.5$, and then $r=n/2+0.25$ and
$r=n/2+0.75$. Although the convergence is slower in this case, the
critical value of $\eta$ where $\lambda_{2}$ exists for infinite
$\mathcal{H}$ is still at $\eta=0$.
\subsection{The third isolated eigenvalue}

Finally, $\lambda_{3}$ exists in a situation that is quite
different from it of $\lambda_{1}$ and $\lambda_{2}$. The pictures
shown below graph the eigenvector denoted as $V_{x}^{(3)}$ with
respect to $\lambda_{3}$ versus the components denoted the index
as $x$ of $V_{x}^{(3)}$ when $n=100$. From top to bottom, $r$ is
increased from $n/2$ to $n/2+1$ by $0.25$ each time, and from left
to right, $\eta$ is increased from $0$ to $2$ by $0.5$ each time.
\begin{figure}[h]
\includegraphics[width=10cm]{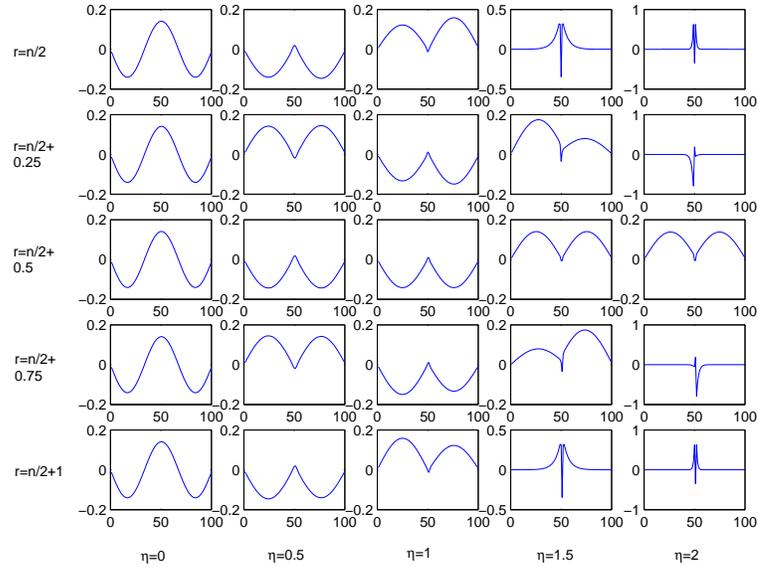}
\centering\caption{$V_{x}^{(3)}$ vs $x$ with n=100}
\end{figure}

As shown on Figure 7 below, for $r=n/2$ and $r=n/2+1$,
$V_{x}^{(3)}$ do not seem to have their concentrated regions or
absolute maxima until $\eta
> 1$, and for $r=n/2+0.25$ and $r=n/2+0.75$, the absolute maxima
appear when $\eta > 1.5$. At $r=n/2+0.5$, $V_{x}^{(3)}$ do not
focus around $x=n/2$ when $n=100$, and the maxima become smaller
and smaller when $n$ increases. Thus, $\lambda_{3}$ may not be an
eigenvalue when $r=n/2+0.5$ because $V_{x}^{(3)}$ may be zero at
$r=n/2+0.5$ for infinite $\mathcal{H}$.

The two tables presented below show the maxima of the square
normalized components in $V_{x}^{(3)}$ at $n=200$ and $n=300$, and
$\eta \in [0.2]$ where $\eta$ is increased by $(2-0)/4=0.5$ every
time, and $r$ ranges from $n/2$ to $n/2+1$ and is increased by
$0.1$ each time.

\begin{center}
\begin{tabular}{|c|c|c|c|c|c|}
  \hline
   & $\eta=0$ & $\eta=0.5$ & $\eta=1$ & $\eta=1.5$ & $\eta=2$ \\
  \hline
   $r=n/2$ & 0.012558 & 0.012896 & 0.119863 & 0.306798 & 0.389453 \\
   $r=n/2+0.1$ & 0.011965 & 0.011392 & 0.116918 & 0.548307 & 0.767640 \\
   $r=n/2+0.2$ & 0.011402 & 0.010295 & 0.014707 & 0.476112 & 0.706397 \\
   $r=n/2+0.3$ & 0.010892 & 0.010118 & 0.015720 & 0.195049 & 0.508429 \\
   $r=n/2+0.4$ & 0.010433 & 0.010182 & 0.013939 & 0.018121 & 0.042290 \\
   $r=n/2+0.5$ & 0.010005 & 0.009957 & 0.009903 & 0.009841 & 0.009763 \\
   $r=n/2+0.6$ & 0.010433 & 0.010182 & 0.013939 & 0.018121 & 0.042290 \\
   $r=n/2+0.7$ & 0.010892 & 0.010118 & 0.015720 & 0.195049 & 0.508429 \\
   $r=n/2+0.8$ & 0.011402 & 0.010295 & 0.014707 & 0.476112 & 0.706397 \\
   $r=n/2+0.9$ & 0.011965 & 0.011392 & 0.116918 & 0.548307 & 0.767640 \\
   $r=n/2+1$ & 0.012558 & 0.012896 & 0.119863 & 0.306798 & 0.389453 \\
  \hline
\end{tabular}
\end{center}

\begin{center}
\begin{tabular}{|c|c|c|c|c|c|}
  \hline
   & $\eta=0$ & $\eta=0.5$ & $\eta=1$ & $\eta=1.5$ & $\eta=2$ \\
  \hline
  $r=n/2$ & 0.008375 & 0.008632 & 0.119863 & 0.306798 & 0.389453 \\
  $r=n/2+0.1$ & 0.007978 & 0.007622 & 0.116918 & 0.548307 & 0.767640 \\
  $r=n/2+0.2$ & 0.007602 & 0.006883 & 0.010120 & 0.476112 & 0.706397 \\
  $r=n/2+0.3$ & 0.007261 & 0.006757 & 0.010572 & 0.195049 & 0.508429 \\
  $r=n/2+0.4$ & 0.006955 & 0.006798 & 0.009331 & 0.012241 & 0.039502 \\
  $r=n/2+0.5$ & 0.006669 & 0.006648 & 0.006624 & 0.006596 & 0.006560 \\
  $r=n/2+0.6$ & 0.006955 & 0.006798 & 0.009331 & 0.012241 & 0.039502 \\
  $r=n/2+0.7$ & 0.007261 & 0.006757 & 0.010572 & 0.195049 & 0.508429 \\
  $r=n/2+0.8$ & 0.007602 & 0.006883 & 0.010120 & 0.476112 & 0.706397 \\
  $r=n/2+0.9$ & 0.007978 & 0.007622 & 0.116918 & 0.548307 & 0.767640 \\
  $r=n/2+1$ & 0.008375 & 0.008632 & 0.119863 & 0.306798 & 0.389453 \\
  \hline
\end{tabular}
\end{center}

\vspace{0.25 in} As presented on the two tables above, at $r=n/2$,
$r=n/2+0.1$, $r=n/2+0.9$, and $r=n/2+1$, the maxima of the square
normalized components in $V_{x}^{(3)}$ do not decrease when $\eta
\geq 1$. Moreover, at $r=n/2+0.2$, $r=n/2+0.3$, $r=n/2+0.7$, and
$r=n/2+0.8$, the maxima do not decrease when $\eta \geq 1.5$.
Lastly, the maxima decrease at $r=n/2+0.4$, $r=n/2+0.5$, and
$r=n/2+0.6$ for all $\eta \in [0,2]$. Based on these two tables,
we may expect that $\lambda_{3}$ exists when $\eta > 1$ while the
interval of $r$ where $\lambda_{3}$ does not exist is becoming
smaller and smaller as $\eta$ is becoming larger and larger.
However, up to now there is not enough argument to support the
claim that $\lambda_{3}$ does not exist at $r=n/2+0.5$. The
pictures below test this claim.

\begin{figure}[h]
\includegraphics[width=10cm]{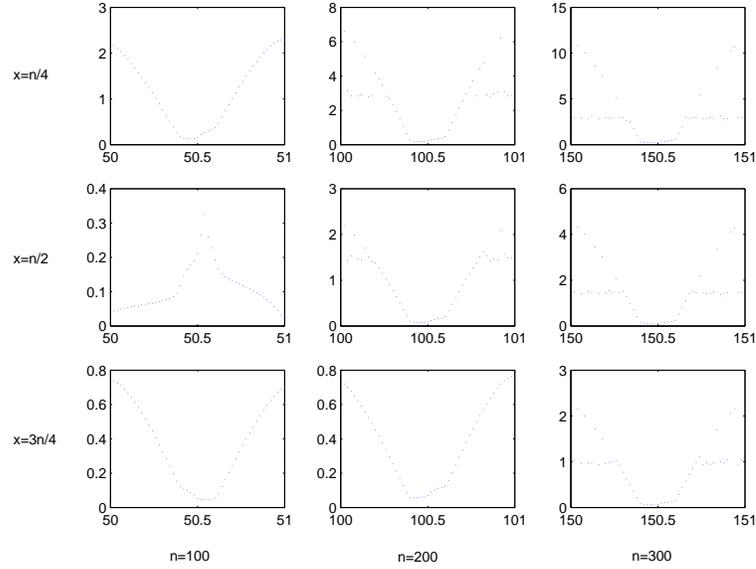}
\centering\caption{$a$ vs $r$ with $\eta=2$}
\end{figure}

Figure 8 below shows the decreasing rates mentioned above at
$x=n/4$, $x=n/2$, and $x=3n/4$ of $V_{x}^{(3)}$ with $n=100$,
$n=200$, and $n=300$. As you can see that the decreasing rate $a$
approaches zero as $r$ tends to $n/2+0.5$. When the decreasing
rate is zero, it implies that the eigenvector at that component
flattens out which does not increase or decrease. Therefore,
$V_{x}^{(3)}$ is expected be zero at $r=n/2+0.5$ when $n$ goes to
infinity.

The graphs of the eigenvectors gives a conclusion that
$\lambda_{3}$ appears when $\eta>1$ while the tables of maxima
shows that it happens when $\eta \geq 1$. So from my numerical
results, I conclude that the critical value of $\eta$ where
$\lambda_{3}$ exists at the borders of $r \in [n/2,n/2+1]$ is
$\eta=1$, and the interval of $r$ where $\lambda_{3}$ does not
exist becomes smaller and smaller as $\eta$ becomes larger and
larger and does not exist at $r=n/2+0.5$. The critical value of
$\eta$ means that when $\eta<1$, $\lambda_{3}$ does not exist for
any $r$, but $\lambda_{3}$ starts to exist when $\eta>1$ for
certain interval of $r$.

\newpage
\section{Conclusion}
My numerical results predict that the spectrum of infinite
dimensional $\mathcal{H}$ contains not only the continuous
spectrum but also the point spectrum. The residual spectrum is
empty because $\mathcal{H}$ is a bounded and self-adjoint
operator. The point spectrum of infinite dimensional $\mathcal{H}$
consists of three eigenvalues with particular regions of $r$ and
$\eta$. The results come from using finite dimensional
$\mathcal{H}$ to approximate infinite dimensional $\mathcal{H}$.
The reason why we can use finite $\mathcal{H}$ to estimate
infinite $\mathcal{H}$ is based on the compact perturbation
theorem because we are already aware of the case of the discrete
Laplacian operator which is also an unperturbed case of
$\mathcal{H}$. Since for finite dimensional $\mathcal{H}$ the
spectrum contains only the point spectrum or a set of eigenvalues,
I focus on the three eigenvalues that are isolated from the set of
the other eigenvalues because those three isolated eigenvalues are
expect to remain in the point spectrum for infinite dimensional
$\mathcal{H}$. These eigenvalues still exist for infinite
$\mathcal{H}$ when their corresponding eigenvectors are nonzero
depending on $r$ and $\eta$. It is sufficient to start with
$r\in[n/2,n/2+1]$ and $\eta \geq 0$ because of the properties of
the diagonal term $\epsilon_{x}$ in $\mathcal{H}$. Based on the
numerics produced by my programs, I conclude that the first
($\lambda_{1}$) and second ($\lambda_{2}$) eigenvalues exist when
$\eta>0$ and for all $r\in[n/2,n/2+1]$. The third eigenvalue
($\lambda_{3}$) exists when $\eta>1$ and the interval of $r$ where
$\lambda_{3}$ does not exist becomes smaller and smaller as $\eta$
becomes larger and larger. Moreover, $\lambda_{3}$ does not exist
at $r=n/2+0.5$.

Some further discussions are listed as follows: First, larger
matrix size $n$ and parameter $\eta$ are needed for a more
complete solution of the existence of $\lambda_{3}$ because my
programs can only conduct the calculations up to $n=300$ and
$\eta=2$. Second, Michoel and Nachtergaele [MN] have proved that
$V_{x}^{(1)}=\frac{1}{\cosh (\eta (x-r))}$. But for $V_{x}^{(2)}$
and $V_{x}^{(3)}$, eigenfunctions like $V_{x}^{(1)}$ may or may
not exist. Third, a detailed analytic proof for using finite
dimensional $\mathcal{H}$ to approximate infinite dimensional
$\mathcal{H}$ is needed because in my paper I state that this path
is possible because of our knowledge of the discrete Laplacian
operator. All these questions are very interested to be studied in
future research.
\newpage

\appendix

\section{Definitions}

\subsection{The definition of a Hilbert space}

\begin{definition}
An inner product on a complex linear space $\mathbf{X}$ is a map
\begin{equation*}
(\cdot , \cdot):\mathbf{X} \times \mathbf{X} \rightarrow
\mathbb{C}
\end{equation*}
such that, for all $x,y,z \in \mathbf{X}$ and $\lambda, \mu \in
\mathbb{C}$:
\begin{enumerate}
    \item $(x,\lambda y + \mu z)=\lambda (x,y) + \mu (x,z)$ (linear in the second
    argument);
    \item $(y,x)=\overline{(x,y)}$ (Hermitian symmetric);
    \item $(x,x)\geq 0$ (nonnegative);
    \item $(x,x)=0$ if and only if $x=0$ (positive definite);
\end{enumerate}
We call a linear space with an inner product an inner product
space or a pre-Hilbert space.
\end{definition}

\begin{definition}
A Hilbert space is a complete inner product space.
\end{definition}

\subsection{The definition of $\ell^{p}(\mathbb{Z})$ and $L^{p}(\mathbb{R}^{n})$}

\begin{definition}
For $1 \leq p < \infty$, the sequence space $\ell^{p}(\mathbb{Z})$
consists of all infinite sequences
$x=(x_{n})_{n=-\infty}^{\infty}$ such that
\begin{equation*}
\sum_{n=-\infty}^{\infty}|x_{n}|^{p} < \infty
\end{equation*}
\end{definition}

\begin{definition}
Suppose that $1 \leq p < \infty$. We denote by $L^{p}(\mathbb{R})$
the set of Lebesgue measurable functions $f:\mathbb{R} \rightarrow
\mathbb{R}$ (or $\mathbb{C}$) such that
\begin{equation*}
\int_{-\infty}^{\infty} |f(x)|^{p} dx < \infty
\end{equation*}
where the integral is a Lebesgue integral, and we identify
functions that differ on a set of measure zero.
\end{definition}

\subsection{The definition of linear, bounded, compact, and self-adjoint operators}

\begin{definition}
A linear map or linear operator $\mathcal{T}$ between real (or
complex) linear spaces $\mathbf{X}$, $\mathbf{Y}$ is a function
$\mathcal{T}:\mathbf{X} \rightarrow \mathbf{Y}$ such that
\begin{equation*}
\mathcal{T}(\lambda x + \mu y) = \lambda \mathcal{T} x + \mu
\mathcal{T} y
\end{equation*}
for all $\lambda$ , $\mu \in \mathbb{R}$ (or $\mathbb{C}$) and
$x$, $y \in \mathbf{X}$. A linear map $\mathcal{T}:\mathbf{X}
\rightarrow \mathbf{Y}$ is called a linear transformation of
$\mathbf{X}$, or a linear operator on $\mathbf{X}$.
\end{definition}

\begin{definition}
Let $\mathbf{X}$ and $\mathbf{Y}$ be two normed linear spaces. We
denote both the $\mathbf{X}$ and $\mathbf{Y}$ norms by $\| \cdot
\|$. A linear map $\mathcal{T}:\mathbf{X} \rightarrow \mathbf{Y}$
is bounded if there is a constant $M \geq 0$ such that
\begin{equation*}
\|\mathcal{T}x\|\leq M\|x\|
\end{equation*}
for all $x\in \mathbf{X}$. If no such constant exists, then we say
that $\mathcal{T}$ is unbounded.
\end{definition}

\begin{definition}
A linear operator $\mathcal{T}:\mathbf{X} \rightarrow \mathbf{Y}$
is compact if and only if every bounded sequence $(x_{n})$ in
$\mathbf{X}$ has a subsequence $(x_{n_{k}})$ such that
$(\mathcal{T}x_{n_{k}})$ converges in $\mathbf{Y}$.
\end{definition}

\begin{definition}
A bounded linear operator
$\mathcal{A}:\mathbf{H}\rightarrow\mathbf{H}$ on a Hilbert space
$\mathbf{H}$ is self-adjoint if and only if
\begin{equation*}
\langle x, \mathcal{A}y \rangle = \langle \mathcal{A}x, y \rangle
\end{equation*}
for all $x,y \in \mathbf{H}$
\end{definition}

\subsection{The definition of kernel and range of an operator}

\begin{definition}
Let $\mathcal{T}:\mathbf{X} \rightarrow \mathbf{Y}$ be a linear
map between linear spaces $\mathbf{X}, \mathbf{Y}$. The null space
or kernel of $\mathcal{T}$, denoted by $ker\mathcal{T}$, is the
subset of $\mathbf{X}$ defined by
\begin{equation*}
ker\mathcal{T}=\{ x\in \mathbf{X} | \mathcal{T}x = 0 \}.
\end{equation*}
The range of $\mathcal{T}$, denoted by $ran\mathcal{T}$, is the
subset of $\mathbf{Y}$ defined by
\begin{equation*}
ran\mathcal{T}=\{ y\in \mathbf{Y} | \exists x\in \mathbf{X} |
\mathcal{T}x=y \}.
\end{equation*}
\end{definition}

\subsection{The definition of the Fourier transform and inverse Fourier transform}

\begin{definition}
The periodic Fourier transform $\mathcal{F}: L^{2}(\mathbb{T})
\rightarrow \ell^{2}(\mathbb{Z})$ that maps a function to its
sequence of Fourier coefficients is defined by
\begin{equation*}
\mathcal{F}f=(\widehat{f}_{n})_{n=-\infty}^{\infty}=\frac{1}{\sqrt{2\pi}}\int_{\mathbb{T}}f(x)e^{-inx}dx
\end{equation*}
such that
\begin{equation*}
f(x)=\frac{1}{\sqrt{2\pi}}\sum_{n=-\infty}^{\infty}\widehat{f_{n}}e^{inx}
\end{equation*}
\end{definition}
\noindent where $\mathbb{T}=[0,2\pi]$.

\begin{definition}
The inverse Fourier transform $\mathcal{F^{*}}:
\ell^{2}(\mathbb{Z}) \rightarrow L^{2}(\mathbb{T})$ is defined by
\begin{equation*}
\mathcal{F^{*}}\widehat{f_{n}}=\frac{1}{2\pi}\sum_{n=-\infty}^{\infty}\widehat{f_{n}}e^{inx}
\end{equation*}
\end{definition}

\subsection{The definition of unitary operators}

\begin{definition}
A linear map $\mathcal{U}: \mathbf{H_{1}} \rightarrow
\mathbf{H_{2}}$ between real or complex Hilbert spaces
$\mathbf{H_{1}}$ and $\mathbf{H_{2}}$ is said to be orthogonal or
unitary, respectively, if it is invertible and if
\begin{equation*}
\langle \mathcal{U}x, \mathcal{U}y \rangle_{\mathbf{H_{2}}} =
\langle x, y \rangle_{\mathbf{H_{1}}}
\end{equation*}
for all $x,y\in \mathbf{H_{1}}$. Two Hilbert spaces
$\mathbf{H_{1}}$ and $\mathbf{H_{2}}$ are isomorphic as Hilbert
spaces if there is a unitary linear map between them.
\end{definition}

\newpage

\section{Programs}
The codes of my three programs: "spectrum.m", "mainresults.m", and
"decreasingrates.m" are presented below. "spectrum.m" is in the
first box, "mainresutls.m" is in the second box and third box, and
"decreasingrates.m" is in the fourth box. The percentage signs
appeared in the codes represent comment for descriptions and not
functioning for certain tasks because "spectgrum.m" plots
$\lambda$ versus $r$ and $\eta$ separately and "mainresults.m"
graphs $V_{x}^{(1)}$, $V_{x}^{(2)}$, $V_{x}^{(3)}$ and calculates
the maxima of the square normalized components of those
eigenvectors at different time as well.

\begin{figure}
\includegraphics{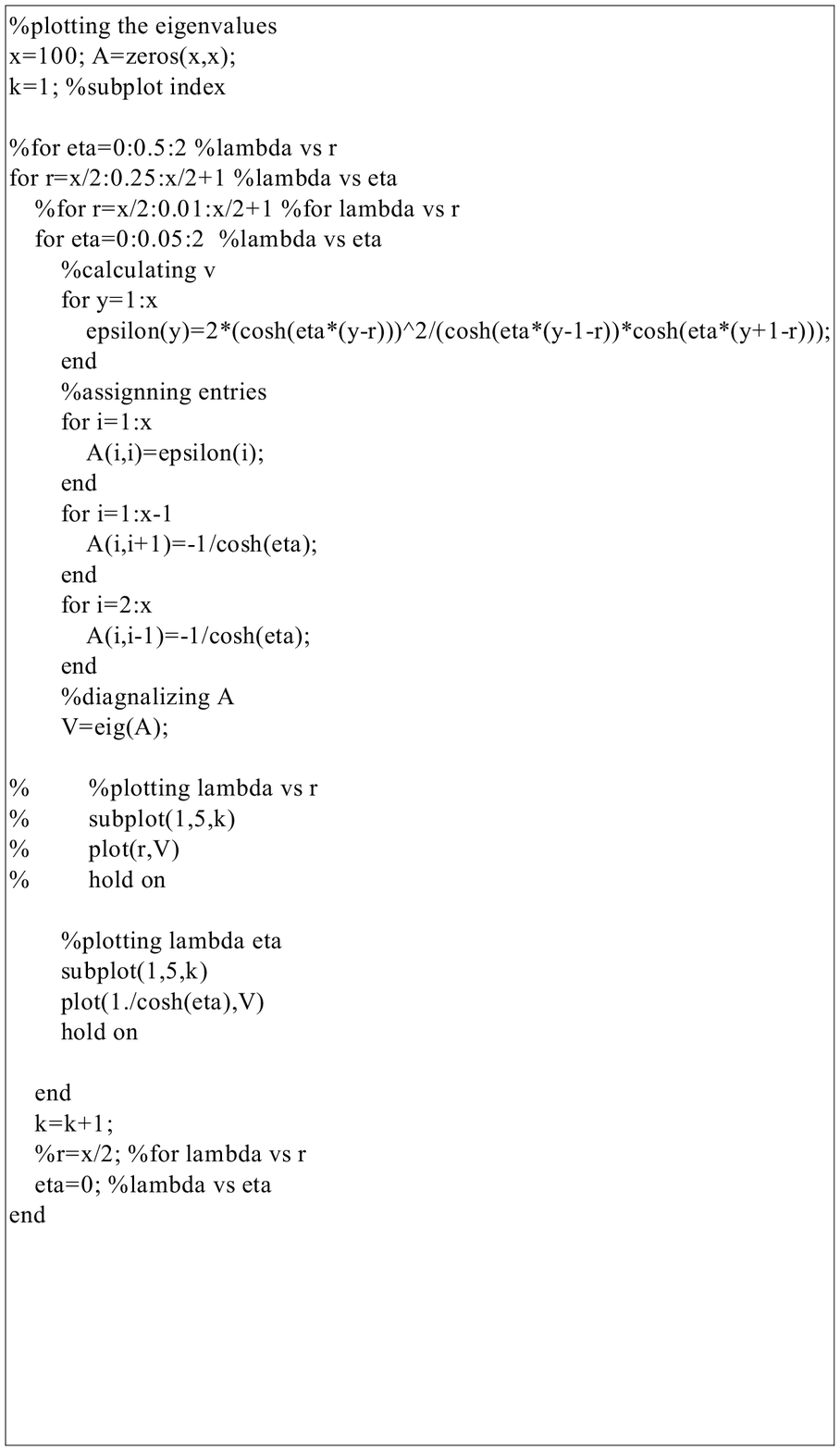}
\centering\caption{spectrum.m}
\end{figure}

\begin{figure}
\includegraphics{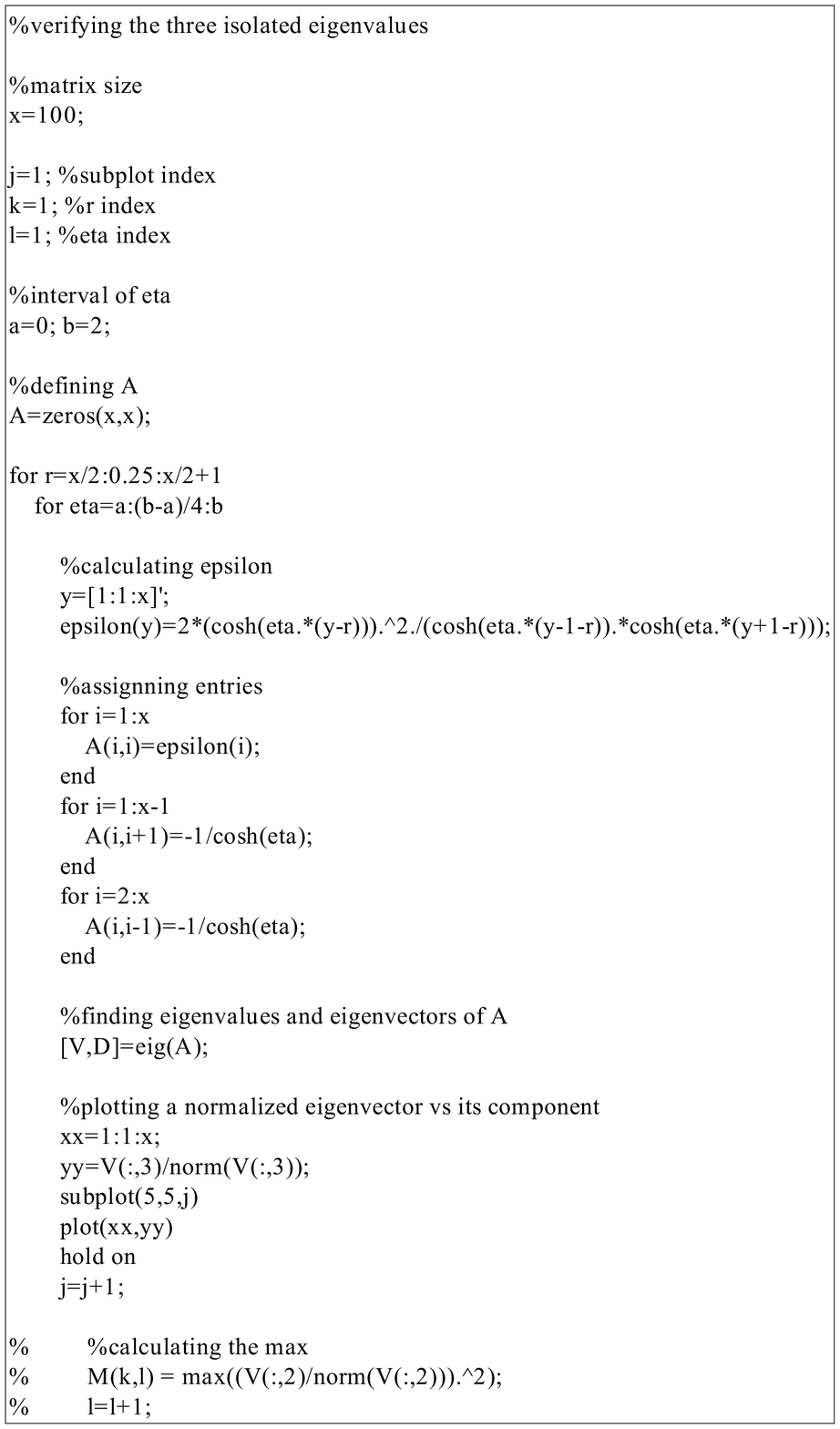}
\centering\caption{mainresults.m}
\end{figure}

\begin{figure}
\includegraphics{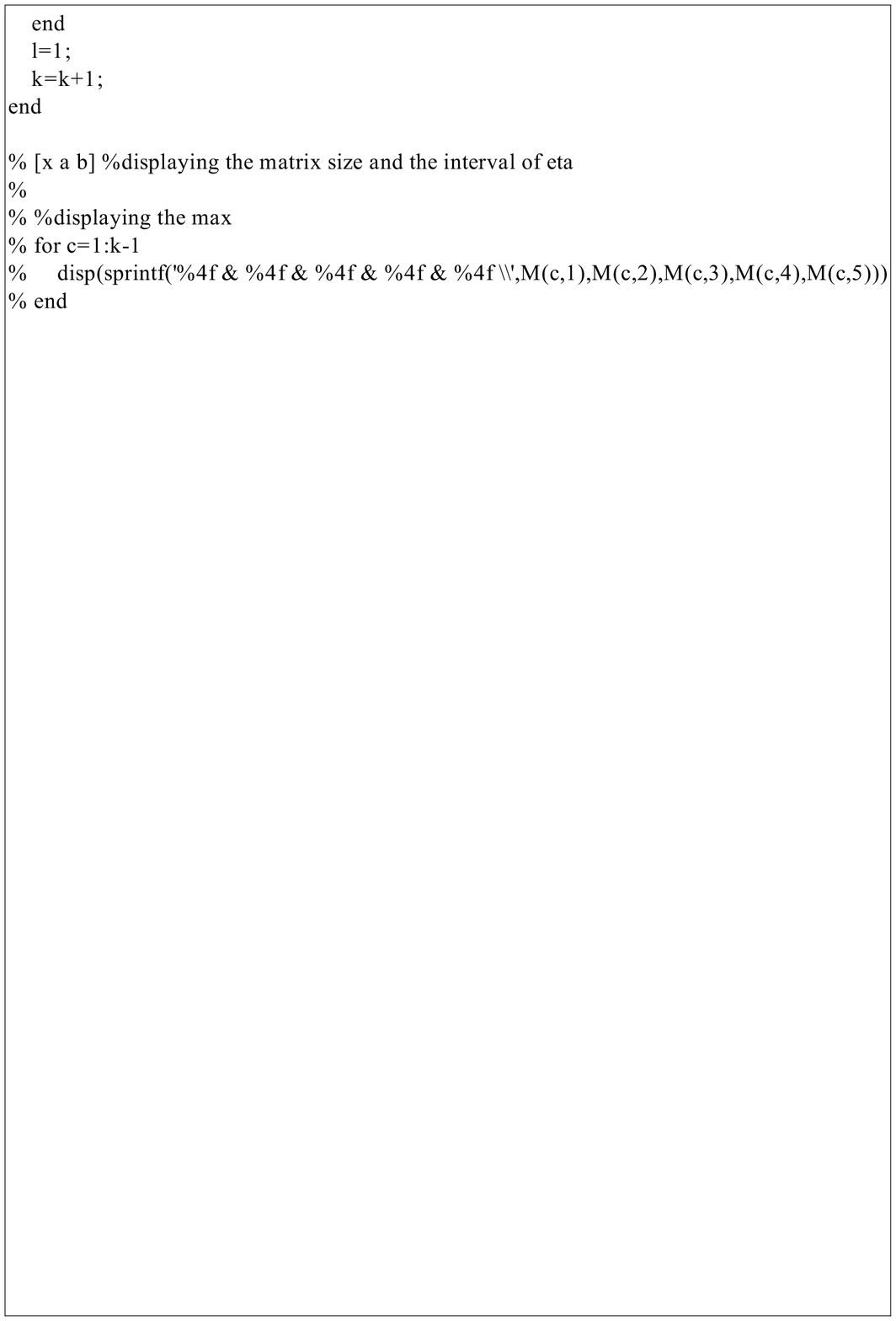}
\centering\caption{mainresults.m cont.}
\end{figure}

\begin{figure}
\includegraphics{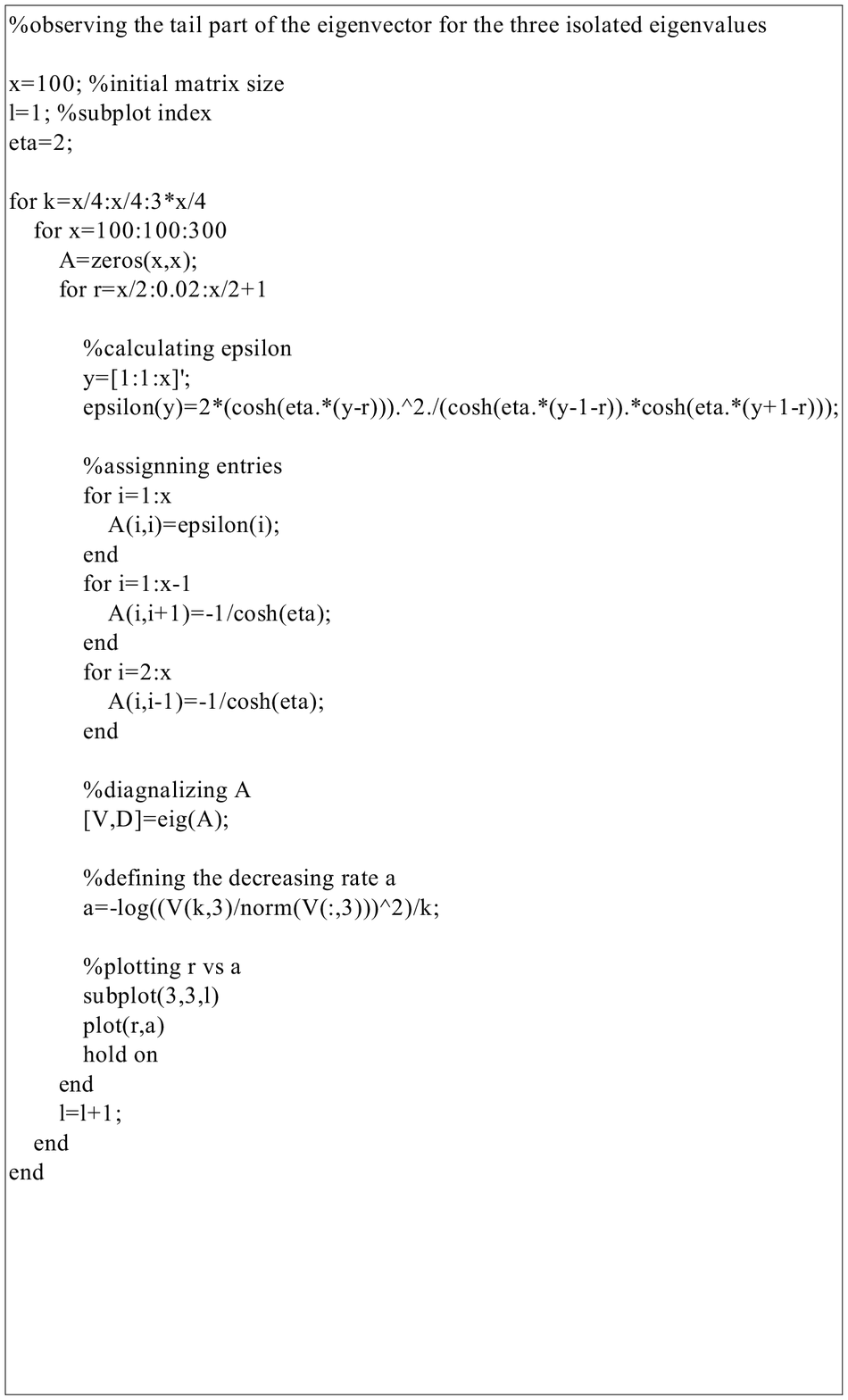}
\centering\caption{decreasingrates.m}
\end{figure}

\newpage


\end{document}